\documentclass[11pt,english]{amsart}
\usepackage{lmodern}
\usepackage{helvet}

\usepackage[T1]{fontenc}
\usepackage[latin9]{inputenc}
\usepackage{amsthm}
\usepackage{amstext}
\usepackage{amssymb}
\usepackage{graphicx}
\usepackage{setspace}
\usepackage{esint}
\makeatletter

\newcommand{\noun}[1]{\textsc{#1}}

\numberwithin{equation}{section}
\numberwithin{figure}{section}
\theoremstyle{plain}
\newtheorem{thm}{\protect\theoremname}[section]
  \theoremstyle{definition}
  \newtheorem{defn}[thm]{\protect\definitionname}
  \theoremstyle{definition}
  \newtheorem{example}[thm]{\protect\examplename}
  \theoremstyle{remark}
  \newtheorem{rem}[thm]{\protect\remarkname}
  \theoremstyle{definition}
  \newtheorem{problem}[thm]{\protect\problemname}
  \theoremstyle{plain}
  \newtheorem{prop}[thm]{\protect\propositionname}

\usepackage{amsfonts}
\usepackage{amssymb}
\usepackage[all]{xy}
\usepackage{array}
\usepackage{lmodern}
\usepackage[T1]{fontenc}

\def\QQ{\mathbb{Q}}
\def\RR{\mathbb{R}}
\def\CC{\mathbb{C}}

\def\ZZ{\mathbb{Z}}
\def\PP{\mathbb{P}}

\def\k12{\mathcal{K}_{\lambda_1,\lambda_2}}
\def\tk12{\tilde{\mathcal{K}}_{\lambda_1,\lambda_2}}
\def\ck12{\check{\mathcal{K}}_{\lambda_1,\lambda_2}}
\def\Res{\text{Res}}

\theoremstyle{definition}

\theoremstyle{definition}

\theoremstyle{theorem}

\theoremstyle{theorem}

\theoremstyle{theorem}

\theoremstyle{definition}

\makeatother

\usepackage{babel}
  \providecommand{\definitionname}{Definition}
  \providecommand{\examplename}{Example}
  \providecommand{\problemname}{Problem}
  \providecommand{\propositionname}{Proposition}
  \providecommand{\remarkname}{Remark}
\providecommand{\theoremname}{Theorem}

\begin{document}

\title{Algebraic cycles and local quantum cohomology}

\author{Charles F. Doran and Matt Kerr }

\begin{abstract}
We review the Hodge theory of some classic examples from mirror symmetry,
with an emphasis on what is intrinsic to the A-model. In particular,
we illustrate the construction of a quantum $\ZZ$-local system on
the cohomology of $K_{\PP^{2}}$ and suggest how this should be related
to the higher algebraic cycles studied in \cite{DK}.
\end{abstract}
\maketitle
This note concerns three types of polarized variations of mixed Hodge
structure (PVMHS) which arise in mirror symmetry:\[\includegraphics[scale=0.6]{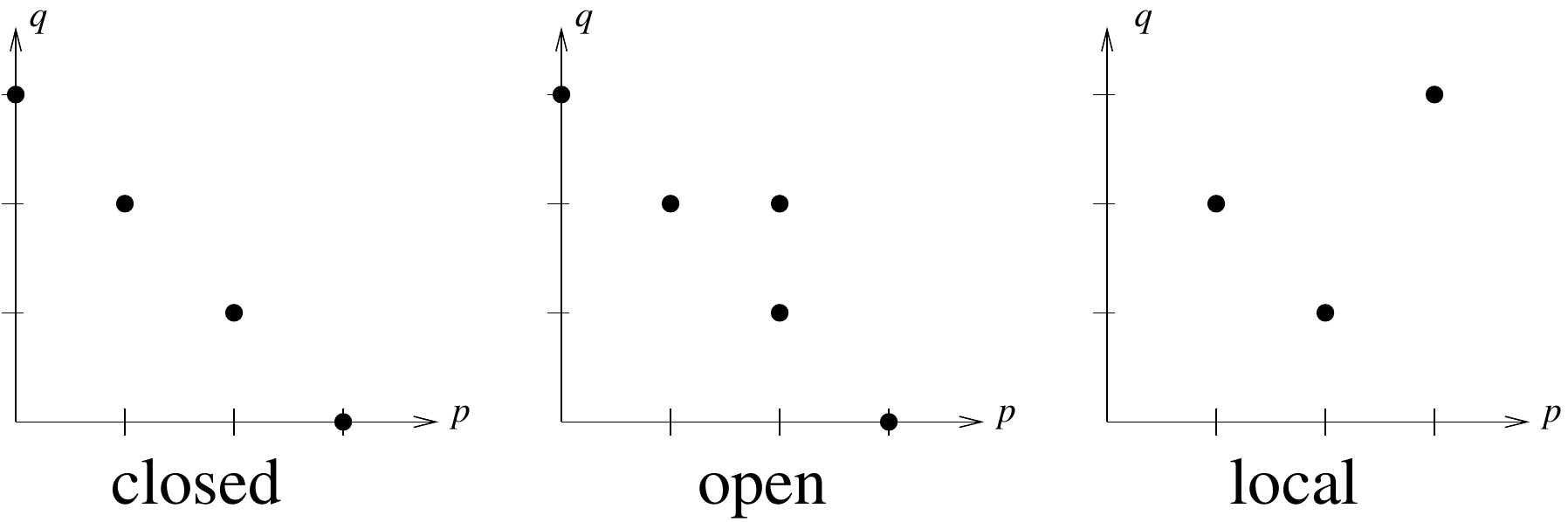}\]In
each case, at the large complex structure boundary point one obtains
a limiting mixed Hodge structure (LMHS) of Hodge-Tate type. It follows
that replacing $W_{\bullet}$ by the relative weight filtration $M_{\bullet}$
produces a new PVMHS of the form\[\includegraphics[scale=0.6]{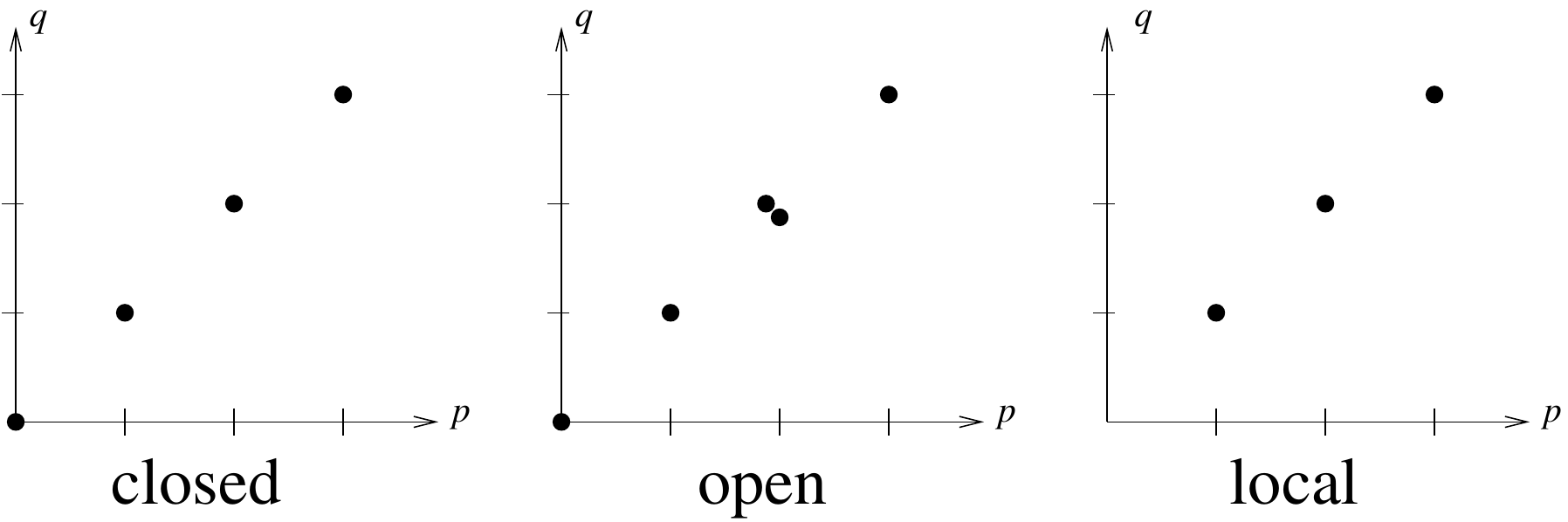}\]occurring
simultaneously in the A and B models. In particular, the $F^{p}\cap M_{2p}$
subspaces identify with $H^{3-p,3-p}$ in quantum cohomology.

Let $\Delta^{*}$ denote the punctured unit disk and write $\mathcal{O}_{\Delta^{*}}=:\mathcal{O}$,
$\Omega_{\Delta^{*}}^{1}=:\Omega^{1}$. A PVMHS $(\mathbb{V},\mathcal{V},\mathcal{F}^{\bullet},W_{\bullet},\nabla,Q)$
over $\Delta^{*}$ comprises
\begin{itemize}
\item a $\ZZ$-local system $\mathbb{V}$ on $\Delta^{*}$,
\item the holomorphic vector bundle $\mathcal{V}$ with sheaf of sections
$\mathbb{V}\otimes\mathcal{O}$,
\item a decreasing filtration by holomorphic subbundles $\mathcal{F}^{j}\subset\mathcal{V}$,
\item an increasing filtration by sub local systems $W_{i}\subset\mathbb{V}_{\QQ}:=\mathbb{V}\otimes\QQ$,
\item a flat connection $\nabla:\mathcal{V}\to\mathcal{V}\otimes\Omega^{1}$
with $\nabla(\mathcal{F}^{\bullet})\subset\mathcal{F}^{\bullet-1}$
and $\nabla(\mathbb{V})=0$, and
\item bilinear forms $Q_{i}:\left(Gr_{i}^{W}\mathbb{V}\right)^{\otimes2}\to\ZZ$,
\end{itemize}
such that each $(Gr_{i}^{W}\mathbb{V}_{s},Gr_{i}^{W}\mathcal{F}_{s}^{\bullet},Q_{i})$
($s\in\Delta^{*}$) yields a polarized Hodge structure. The PVMHS
considered here, as well as all PVMHS arising from geometry, are \emph{admissible}
-- i.e. have well-defined LMHS at $0$. 

In the above pictures, the number of bullets in position $(p,q)$
signifies the dimension of the summand in the Deligne bigrading on
$\mathcal{V}$ defined pointwise by
\[
I^{p,q}(\mathcal{V}_{s}):=F^{p}\cap W_{p+q}\cap\left(\overline{F^{q}}+\sum_{j\geq0}\left\{ \overline{F^{q-j-1}}\cap W_{p+q-j-2}\right\} \right).
\]
This bigrading is uniquely determined by the properties
\begin{enumerate}
\item $\oplus_{p\geq j}\oplus_{q}I^{p,q}(\mathcal{V}_{s})=\mathcal{F}_{s}^{j}$
\item $\oplus_{p+q\leq i}I^{p,q}(\mathcal{V}_{s})=(W_{i})_{s}\otimes\CC$
\item $\overline{I^{b,a}(\mathcal{V}_{s})}\equiv I^{a,b}(\mathcal{V}_{s})$
modulo $\oplus_{p<a}\oplus_{q<b}I^{p,q}(\mathcal{V}_{s})$.
\end{enumerate}
In passing to the limit, heuristically one may visualize the bullets
in each line $p+q=i$ moving up and down in such a way that the end
result remains symmetric about this line.

\subsubsection*{Notation:}

Set $\ell(s):=\frac{\log(s)}{2\pi i}$. We shall often write $\mathcal{V}$
(instead of the 6-tuple) for a PVMHS.

\subsubsection*{Acknowledgments:}

We thank E. Zaslow for a helpful conversation. The first author wishes
to recognize support from the NSERC Discovery Grants Program and the
second author from the NSF under Standard Grant DMS-1068974.

\section{Closed string}

Beginning on the B-model side, recall how the LMHS construction works
for a pure ($\mathbb{Z}$-)VHS $\mathcal{V}$ of weight 3 over $\Delta^{*}$
with unimodular polarization $Q$. The weight filtration is the trivial
one $W_{3}=\mathcal{V}\supset W_{2}=\{0\}$. Denote the (unipotent
part of the) monodromy operator by $T$, with nilpotent logarithm
\[
N:=\log(T):\,\mathbb{V_{Q}}\to\mathbb{V_{Q}}.
\]
There exists an unique filtration
\[
M_{-1}=\{0\}\subset M_{0}\subset M_{1}\subset\cdots\subset M_{6}=\mathbb{V_{Q}}
\]
satisfying $N(M_{\alpha})\subset M_{\alpha-2}$ and $N^{\ell}:Gr_{3+\ell}^{M}\overset{\cong}{\to}Gr_{3-\ell}^{M}$.
Untwisting the local system by
\[
\tilde{\mathbb{V}}:=e^{-\ell(s)N}\mathbb{V},
\]
we obtain the canonical extension 
\[
\mathcal{V}_{e}:=\tilde{\mathbb{V}}\otimes\mathcal{O}_{\Delta}.
\]
Let $\{\gamma_{i}\}$ be a multivalued basis of $\mathbb{V}$ generating
the steps of the integral filtration $M_{m}^{\mathbb{Z}}:=\mathbb{V}\cap M_{m}$,
and set $\tilde{\gamma}_{i}:=e^{-\ell(s)N}\gamma_{i}\in\Gamma(\Delta,\tilde{\mathbb{V}})$.
\begin{defn}
The LMHS of $\mathcal{V}$, denoted informally $V_{lim}$, is given
by the data $V_{lim}^{\mathbb{Z}}:=\mathbb{Z}\langle\{\tilde{\gamma}_{i}(0)\}\rangle$,
$\mathcal{F}_{lim}^{\bullet}:=\mathcal{F}_{e}^{p}(0)$, and (monodromy
weight) filtration $M_{\bullet}$ on $V_{lim}:=\mathcal{V}_{e}(0)$.
\end{defn}
Assume that $V_{lim}$ is \emph{Hodge-Tate}, i.e. $Gr_{2j}^{M}\cong\mathbb{Z}(-j)^{\oplus d_{j}}$
for $j=0,1,2,3$ and $\{0\}$ otherwise. (For example, the LMHS for
$H^{3}$ of the quintic mirror is of this type, while that for the
Fermat quintic family is \emph{not}.) In the rank $4$ setting, where
we must have all $d_{j}=1$, we may pick (for each $j$) a holomorphic
section $e_{j}\in\Gamma(\Delta,\mathcal{F}_{e}^{j}\cap M_{2j}^{\mathbb{C}})$
mapping to the image of $\gamma_{j}\in\Gamma(\mathfrak{H},M_{2j}^{\mathbb{Z}})$
in $\Gamma(\Delta^{*},Gr_{2j}^{M}\mathcal{V})$ hence generating the
latter. Write $e=\{e_{3},e_{2},e_{1},e_{0}\}$ and $\gamma=\{\gamma_{3},\gamma_{2},\gamma_{1},\gamma_{0}\}$
for the two bases.

To make things explicit, we have (for some $a,b\in\mathbb{Z}$ and
$e,f\in\mathbb{Q}$)\begin{equation}\label{Q and N}[Q]_{\gamma}=\left(\begin{array}{cccc}0 & 0 & 0 & 1\\0 & 0 & 1 & 0\\0 & -1 & 0 & 0\\-1 & 0 & 0 & 0\end{array}\right)=[Q]_{e}\;\text{ and }\;[N]_{\gamma}=\left(\begin{array}{cccc}0 & 0 & 0 & 0\\a & 0 & 0 & 0\\e & b & 0 & 0\\f & e & -a & 0\end{array}\right)
\end{equation} (cf. \cite{GGK}), in which we shall demand that $|a|=1$. Replacing
the local coordinate $s$ by $q:=e^{2\pi\sqrt{-1}\tau}$, where $\tau:=Q(\gamma_{1},e_{3})$,
and making full use of the bilinear relations (e.g. $Q(\mathcal{F}^{1},\mathcal{F}^{3})=0=Q(\mathcal{F}^{2},\mathcal{F}^{2})$),
the limiting period matrix becomes (cf. {[}op. cit.{]}) \begin{equation}\label{LMHS}_{\tilde{\gamma}(0)}[\mathbf{1}]_{e(0)}=\left(\begin{array}{cccc}1 & 0 & 0 & 0\\0 & 1 & 0 & 0\\\frac{f}{2} & e & 1 & 0\\\alpha_{0} & \frac{f}{2} & 0 & 1\end{array}\right).
\end{equation}
\begin{example}
\label{mir quintic constants}For the mirror quintic family, we have
(cf. {[}op. cit.{]}, where the computation is based on \cite{CdOGP})
$a=-1$, $b=5$, $e=\frac{11}{2}$, $f=-\frac{25}{6}$, and $\alpha_{0}=\frac{25i}{\pi^{3}}\zeta(3)=:C$.
\end{example}
Following Deligne \cite{De}, the $e_{j}(q)|_{\Delta^{*}}$ provide
the Hodge(-Tate) basis of a PVMHS $(\mathbb{V},\mathcal{V},\mathcal{F}^{\bullet},M_{\bullet},\nabla)$
on $\Delta^{*}$, denoted $\mathcal{V}_{rel}$ for short. For the
connection, we have
\[
[\nabla]_{e}=d+\left(\begin{array}{cccc}
0 & 0 & 0 & 0\\
1 & 0 & 0 & 0\\
0 & -Y(q) & 0 & 0\\
0 & 0 & -1 & 0
\end{array}\right)\otimes\frac{dq}{(2\pi\sqrt{-1})q}
\]
where $Y(q)$ defines the Yukawa coupling. In the event that $\mathcal{V}$
comes from $H^{3}(X)$, and $\Phi$ denotes the Gromov-Witten prepotential
of the mirror $X^{\circ}$ (composed with the inverse mirror map),
according to mirror symmetry we have $Y=\Phi''':=\frac{d^{3}\Phi}{d\tau^{3}}.$
\begin{example}
\label{exa:The-mirror-quintic}The mirror quintic VHS arises from
$H^{3}$ of $X_{\xi}$, which is a smooth compactification of 
\[
\left\{ 1-\xi\left(\sum_{i=1}^{4}x_{i}+\frac{1}{\prod_{i=1}^{4}x_{i}}\right)=0\right\} \subset\left(\mathbb{C}^{*}\right)^{\times4}.
\]
Taking $s:=\xi^{5}$, we obtain $\tau$ and $q$ as above, and 
\[
\Phi(q)=\frac{5}{6}\tau^{3}+\Phi_{h}(q),
\]
where the holomorphic part 
\[
\Phi_{h}(q)=\frac{1}{(2\pi i)^{3}}\sum_{d\geq1}N_{d}q^{d}.
\]
From \cite{CdOGP,GGK,Pe}, we have the mixed Hodge basis \begin{flalign*}
\; & e_0 = \gamma_0 &\\
\; & e_1 = \gamma_1 - \tau \gamma_0 &\\
\; & e_2 = \gamma_2 - \left( 5\tau+\frac{11}{2}+\Phi_h'' \right) \gamma_1 + \left( \frac{5}{2}\tau^2 + \frac{25}{12} + \tau\Phi_h'' - \Phi_h'\right) \gamma_0 &\\
\; & e_3 = \gamma_3 + \tau \gamma_2 - \left( \frac{5}{2}\tau^2 +\frac{11}{2}\tau -\frac{25}{12}+\Phi_h'\right)\gamma_1 &\\
\; & \mspace{200mu} +\left(\frac{5}{6}\tau^3 + \frac{25}{12}\tau-C+\tau\Phi_h'-2\Phi_h \right)\gamma_0 .
\end{flalign*}Here $e_{3}$ can also be viewed as the class of a holomorphic $3$-form
in the original VHS, whose LMHS is reflected by the presence of $C$.
The mirror $X^{\circ}$ is the Fermat quintic.
\end{example}
Turning to the A-model, we need to define an integral structure, Hodge
and weight filtrations on 
\[
H^{\text{even}}(X^{\circ})=H^{3,3}\oplus H^{2,2}\oplus H^{1,1}\oplus H^{0,0}
\]
which will lead to VHS, LMHS, and VMHS isomorphic to those on $H^{3}(X)$.
These variations will be defined over a small disk $0<|q|<\epsilon$.
For constructing them, the general idea is to use the family of algebraic
structures on $H^{even}$ parametrized by $\tau[H]\in H^{1,1}(X^{\circ})$,
known as the \emph{(small) quantum cohomology}. (Here $[H]$ the the
class of a hyperplane section and $\tau=\ell(q)$, and we are working
in the rank $4$ setting.)

For the filtrations, we set 
\[
F^{a}H^{\text{even}}=\oplus_{i\leq3-a}H^{i,i}\;,\;\;\; M_{b}H^{\text{even}}=\oplus_{j\geq3-\frac{b}{2}}H^{j,j}
\]
so that $\mathcal{F}^{3-k}\cap M_{6-2k}=H^{i,i}(X^{\circ},\mathbb{C})$
as a subspace of $H^{\text{even}}$. This is where the ``naive''
fundamental classes of coherent sheaves or algebraic cycles of codimension
$i$ lie. In contrast, the integral local system will be generated
by \emph{quantum-deformed} fundamental classes of algebraic cycles
on $X^{\circ}$. Alternately, we can regard the flat structure as
given by the solution to a quantum differential equation 
\[
\nabla=d+E\otimes\frac{dq}{(2\pi\sqrt{-1})q},
\]
which gives the integral structure up to a constant. (Note that $d$
differentiates with respect to $\oplus_{i}H^{i,i}(X^{\circ},\mathbb{C})$.)
Since $E$ kills $M$-graded pieces, we get a natural identification
between $Gr_{2i}^{M}$ of this ``integral structure'' and $H^{i,i}(X^{\circ},\mathbb{Z})$.
\begin{example}
\label{integral basis}For $X^{\circ}$ the Fermat quintic, we have
Hodge basis 
\[
[X^{\circ}]=e_{3},\;[H]=e_{2},\;-[L]=e_{1},\;[p]=e_{0}
\]
 where $H$ is a hyperplane section, $L$ a line and $p$ a point.
The minus sign on $[L]$ ensures that the form 
\[
Q(\alpha,\beta):=(-1)^{\frac{\deg(\alpha)}{2}}\int_{X^{\circ}}\alpha\cup\beta
\]
has matrix $[Q]_{e}$ as above, which is necessary for equality of
\emph{polarized} VHS.

For the quantum deformed classes, we invert the relations of Example
\ref{exa:The-mirror-quintic} to obtain \begin{flalign*}
\; & [X^{\circ}]_{\mathcal{Q}}=\gamma_3=[X^{\circ}] - \tau[H]+\left( \frac{5}{2}\tau^2 +\frac{25}{12}+\tau\Phi_h''-\Phi_h' \right)[L] &\\
\;   & \mspace{200mu} + \left( -\frac{5}{6}\tau^3-\frac{25}{12}\tau +C -\tau\Phi_h'+2\Phi_h\right)[p], &\\
\; & [H]_{\mathcal{Q}} = \gamma_2 = [H] - \left( 5\tau + \frac{11}{2}+\Phi_h''\right) [L] + \left( \frac{5}{2}\tau^2 + \frac{11}{2}\tau-\frac{25}{12} + \Phi_h'\right)[p], &\\
\; & [L]_{\mathcal{Q}} = -\gamma_1 = [L] - \tau [p], &\\
\; & [p]_{\mathcal{Q}} = \gamma_0 = [p].
\end{flalign*}These are solutions to the above differential equation with $E$ given
by the (small) quantum product $[H]*$ defined by
\[
[H]*[X^{\circ}]=[H],\;\;[H]*[H]=\Phi'''[L],\;\;[H]*[L]=[p],\text{ and }[H]*[p]=0.
\]
(Note that this is consistent with cup product, in the sense that
$[H]\cup[H]=5[L]=\Phi'''(0)[L]$.) The resulting variations of HS
on $H^{\text{even}}(X^{\circ})$ and $H^{3}(X)$ match by construction.
\end{example}
The natural question at this point is: \emph{how much of this ``common
$\ZZ$-VHS'' is intrinsic to the A-model, and not just the B-model?}
Clearly the issue lies not in the Hodge and monodromy weight filtrations
(given by the grading of $H^{even}$ by degree), or the polarizing
form $Q$, or the $\nabla$-flat complex local system (given by the
quantum product), but in the integral structure on the latter. Another
way to think of this (cf. \cite{De}) is that we must determine the
``constant of integration'' of the VHS, or equivalently the LMHS
\eqref{LMHS}.

Naively, one could try to find a basis $\delta$ of the local system
with integral $[Q]_{\delta}$ and integral monodromy matrices (which
are computable in principle by analytic continuation). Unfortunately
the result may not be unique, even after identifying bases related
by a rational symplectic matrix. In the above example, one could have
\[
\delta_{3}=\frac{\gamma_{3}}{\sqrt{5}}+\frac{\gamma_{2}}{\sqrt{5}}\,,\;\delta_{2}=\frac{\gamma_{2}}{\sqrt{5}}-\frac{3\gamma_{1}}{\sqrt{5}}-\frac{3\gamma_{0}}{\sqrt{5}}\,,\;\delta_{1}=\sqrt{5}\xi_{1}\,,\;\delta_{0}=\sqrt{5}\gamma_{0},
\]
which produces the (distinct) quintic \emph{twin} mirror $\ZZ$-VHS.
Indeed, in \cite{DM} this phenomenon is responsible for the bifurcation
of each $\RR$-VHS into finitely many distinct $\ZZ$-VHS.

Instead, what is needed is a direct construction of an integral structure
on quantum cohomology, which has only recently been realized by Iritani
\cite{Ir1,Ir2} and Katzarkov-Kontsevich-Pantev \cite{KKP}. We illustrate
how this works in the setting where $X^{\circ}$ is a smooth CY 3-fold,
and $\dim H^{even}(X^{\circ})=4$. A map $\sigma$ from $H^{even}$
to multivalued $\nabla$-flat sections (in a neighborhood of $q=0$),
defined in terms of Gromov-Witten theory, has been known for some
time (cf. \cite[secs. 8.5.3, 10.2.2]{CK}). If $\alpha_{i}\in H^{2(3-i)}(X^{\circ})$
($i=0,1,2,3$) denote a $Q$-symplectic basis with $\alpha_{2}=[H]$,
this boils down to first setting 
\[
\tilde{\sigma}(\alpha_{0}):=\alpha_{0},\;\;\tilde{\sigma}(\alpha_{1}):=\alpha_{1},\;\;\tilde{\sigma}(\alpha_{2}):=\alpha_{2}+\Phi_{h}''\alpha_{1}+\Phi_{h}'\alpha_{0},
\]
\[
\tilde{\sigma}(\alpha_{3}):=\alpha_{3}+\Phi_{h}'\alpha_{1}+2\Phi_{h}\alpha_{0}
\]
and then 
\[
\sigma(\alpha):=\tilde{\sigma}\left(e^{-\tau[H]}\cup\alpha\right):=\sum_{k\geq0}\frac{(-1)^{k}}{k!}\tilde{\sigma}\left([H]^{k}\cup\alpha\right).
\]
(In our running example, we obviously have in mind $\alpha_{3}=[X^{\circ}]$,
$\alpha_{2}=[H]$, $\alpha_{1}=-[L]$, and $\alpha_{0}=[p]$.) These
are $\nabla$-flat sections with monodromy \begin{equation}\label{monodromy}T(\sigma(\alpha))=\sigma\left(e^{-[H]}\cup\alpha\right).
\end{equation}We also set $\sigma_{\infty}(\alpha):=\tilde{\sigma}(\alpha)|_{q=0}.$

The key new ingredient introduced by \cite{Ir1,KKP} is a characteristic
class defined using the $\Gamma$-function, and which in our setting
specializes to \begin{equation}\label{gamma class}\hat{\Gamma}(X^{\circ}):=\exp\left(\sum_{k\geq2}\frac{(-1)^{k}(k-1)!}{(2\pi i)^{k}}\zeta(k)ch_{k}(TX^{\circ})\right)\in H^{even}(X^{\circ}).
\end{equation} Using it, we may assign a flat section \begin{equation}\label{gamma Z-str}\gamma(\xi):=\sigma\left(\hat{\Gamma}(X^{\circ})\cup ch(\xi)\right)
\end{equation}to each $\xi\in K_{0}^{num}(X^{\circ})$, which defines a $\ZZ$-local
system. (Similarly, we can define $\tilde{\gamma}(\xi)$, $\gamma_{\infty}(\xi)$
by applying $\tilde{\sigma}$, $\sigma_{\infty}$.) A strong indication
that $\hat{\Gamma}$ gives the right ``correction'' is Iritani's
result (cf. \cite[Prop. 2.10]{Ir1}) that the Mukai pairing
\[
\left\langle \xi,\xi'\right\rangle :=\int_{X^{\circ}}ch(\xi^{\vee}\otimes\xi')\cup Td(X^{\circ})\;=\; Q(\gamma(\xi),\gamma(\xi')).
\]
Moreover, since $ch(\mathcal{O}(-1))=e^{-[H]}$, \eqref{monodromy}
implies that 
\[
T(\gamma(\xi))=\gamma(\mathcal{O}(-1)\otimes\xi)
\]
--- an elementary example of how a categorical autoequivalence of
$D^{b}(X^{\circ})$ corresponds to monodromy. The autoequivalences
corresponding to monodromies arising away from $q=0$ have been explicitly
identified in \cite{CIR}.
\begin{example}
Once more we take $X^{\circ}$ to be the Fermat quintic, which has
total Chern class $c(X^{\circ})=1+50[L]-200[p]$ and Todd class $Td(X^{\circ})=1+\frac{25}{6}[L]$.
A Mukai-symplectic basis of $K_{0}^{num}(X^{\circ})$ is
\[
\xi_{3}:=\mathcal{O}_{X^{\circ}},\;\xi_{2}:=\mathcal{O}_{H}-3\mathcal{O}_{L}-8\mathcal{O}_{p},\;\xi_{1}:=-\mathcal{O}_{L}-\mathcal{O}_{p}\equiv-\mathcal{O}_{L}(1),\;\xi_{0}:=\mathcal{O}_{p};
\]
this in fact (referring to Example \ref{mir quintic constants} and
\eqref{Q and N}) satisfies $[\mathcal{O}(-1)\otimes]_{\xi}=\exp\left([N]_{\gamma}\right)$.
(Note that taking $\xi_{2}=\mathcal{O}_{H}$ and $\xi_{1}=\mathcal{O}_{L}$
does \emph{not} yield a symplectic basis.) From 
\[
ch(\xi_{3})=[X^{\circ}],\; ch(\xi_{2})=[H]-\frac{11}{2}[L]-\frac{25}{6}[p],\; ch(\xi_{1})=-[L],\; ch(\xi_{0})=[p]
\]
and $\hat{\Gamma}(X^{\circ})=[X^{\circ}]+\frac{25}{12}[L]+C[p]$,
a straightforward computation gives that
\[
\gamma(\xi_{i})=\gamma_{i}\;\;\;(i=0,1,2,3),
\]
with the $\{\gamma_{i}\}$ exactly as in Example \ref{integral basis}.
Moreover, the $\{\gamma_{\infty}(\xi_{i})\}$ recover the LMHS matrix
\eqref{LMHS} (with $e,f,\alpha_{0}$ as in Example \ref{mir quintic constants}),
including the crucial constant $C$ which visibly comes from $\hat{\Gamma}$.\end{example}
\begin{rem}
\label{rem arithmetic}The toric-hypersurface CY 3-fold families from
which B-model VHS's are often produced are intrinsically defined over
$\QQ$. Moreover, by virtue of its toric nature, the large complex
structure limit may be regarded as a $\QQ$-semistable degeneration.
The general conjectural framework surrounding the limiting motive
(cf. \cite[(III.B.5)]{GGK}) therefore predicts that the class $\alpha_{0}\in Ext_{\text{MHS}}^{1}(\QQ(-3),\QQ(0))\cong\CC/\QQ$
arising in the corresponding LMHS is always a rational multiple of
the constant $\frac{\zeta(3)}{(2\pi i)^{3}}$, motivating its appearance
in \eqref{gamma class}.%
\footnote{Note that we are interested in the arithmetic of locally complete
CY families; taking irrational ``slices'' of such to force an extension
both misses the point and will not affect $\alpha_{0}$.%
}

The ``non-toric'' degenerations at the conifold and Gepner points,
on the other hand, produce singular fibers whose desingularization
may introduce an algebraic extension of $\QQ$, leading to an arithmetically
richer LMHS. One should try to use mirror symmetry to get at this,
perhaps beginning with\end{rem}
\begin{problem}
Adapt the (A-model) $\hat{\Gamma}$-integral structure on FJRW theory
introduced in \cite{CIR} to the explicit computation of the periods
of (B-model) LMHS at the Gepner point ($s=\infty$).
\end{problem}
See $\S4$ for another source of algebraic extensions.

\section{Local string}

This section is based on a simple example studied by \cite{CKYZ},
\cite{MOY}, \cite{Ho}, and \cite{DK}. Once and for all we set\begin{equation}\label{eqn HV mir}Y_{\xi}:=\left\{ (x,y;u,v)\in (\CC^*)^2 \times \CC^2 \right. \left| 1-\xi\left(x+y+\frac{1}{xy}\right)+u^{2}+v^{2}=0\right\} ,
\end{equation}the so-called \emph{Hori-Vafa mirror} of $Y^{\circ}=K_{\mathbb{P}^{2}}$.
The canonical holomorphic $(3,0)$ form on $Y_{\xi}$ is given by
\[
\eta_{\xi}=2\sqrt{-1}\Res_{Y_{\xi}}\left(\frac{\frac{dx}{x}\wedge\frac{dy}{y}\wedge du\wedge dv}{1-\xi(x+y+\frac{1}{xy})+u^{2}+v^{2}}\right).
\]
The $3$-cycles are spanned in homology by (a) a real $3$-torus $\mathbb{T}^{3}$
and (b) circle-bundles over membranes in $(\mathbb{C}^{*})^{\times2}$
bounding $1$-cycles on the thrice-punctured elliptic curve $W_{\xi}^{*}:=\left\{ (x,y)\in (\CC^*)^2 \right. \left| 1-\xi(x+y+\frac{1}{xy})=0\right\}.
$ The circle is pinched to a point over the $1$-cycles. \[\includegraphics[scale=0.5]{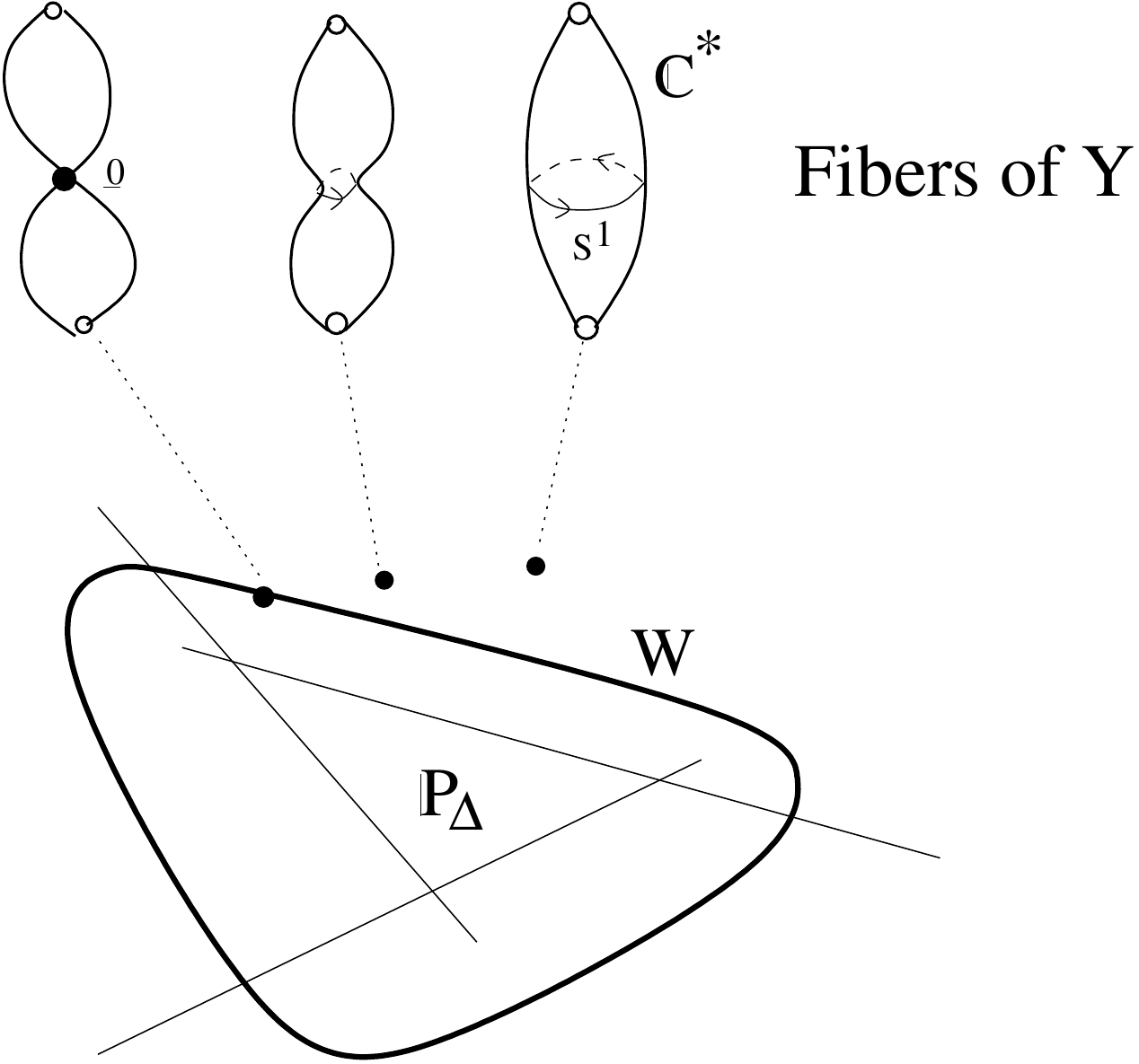}\]We
write $W_{\xi}$ for the complete elliptic curve, $\tilde{\omega}_{\xi}:=\frac{1}{2\pi i}\Res_{W_{\xi}}\left(\frac{\frac{dx}{x}\wedge\frac{dy}{y}}{1-\xi(x+y+\frac{1}{xy})}\right)$
for the canonical holomorphic $1$-form, and $\varphi_{0}$,$\varphi_{1}$
for $1$-cycles spanning $H_{1}(W_{\xi},\mathbb{Z})$ with periods
$\pi_{i}:=\int_{\varphi_{i}}\tilde{\omega}_{\xi}$. In particular,
we let $\varphi_{0}$ be the vanishing cycle and $\omega_{\xi}:=\tilde{\omega}_{\xi}/\pi_{0}$
the normalization of the $1$-form so that $\int_{\varphi_{0}}\omega_{\xi}\equiv1$.

Denoting the membrane construction (b) by $\mathcal{M}$, we have
the short exact sequence \[
\xymatrix{ 0 \ar [r] & \mathbb{Z} \langle \mathbb{T}^3 \rangle \ar [r] & H_3(Y)\ar [r] & \ker\left\{ H_1(W^*) \to H_1((\mathbb{C}^*)^2) \right\}(1) \ar [r] \ar [d]_{\cong} \ar @/^2pc/ @{-->} [l]^{\mathcal{M}} & 0 \\ & & & H_1(W)(1) }
\](cf. {[}DK, sec. 5{]}).%
\footnote{The isomorphism is valid only rationally, but can be made integral
by replacing $H_{1}(W,\ZZ)$ by $\ZZ\langle3\varphi_{0},\varphi_{1}\rangle$,
which is done tacitly below.%
} Its dual\[
\xymatrix{0 & \mathbb{Z}(-3) \ar [l] & H^3(Y) \ar [l] & H^1(W)(-1) \ar [l]_{\mu} & 0 \ar [l]}
\]yields an extension class 
\[
\varepsilon\in\text{Ext}_{\text{MHS}}^{1}\left(\mathbb{Z}(-2),H^{1}(W)\right)\cong\text{Hom}\left(H_{1}(W),\mathbb{C}/\mathbb{Z}(2)\right).
\]
Miraculously, this is the image of a higher cycle $\Xi\in K_{2}^{\text{alg}}(W)$
by a generalized Abel-Jacobi map \cite{DK}, and the periods of $\eta$
may be described by
\[
\frac{1}{2\pi\sqrt{-1}}\int_{\mathcal{M}(\gamma)}\eta\underset{\mathbb{Z}(2)}{\equiv}\langle AJ(\Xi),\gamma\rangle_{W}\;,\;\;\;\;\frac{1}{(2\pi\sqrt{-1})^{3}}\int_{\mathbb{T}^{3}}\eta=1.
\]

Normalizing the local coordinate $s:=\xi^{3}$ to $q$ where 
\[
\ell(q):=\tau:=\frac{\pi_{1}}{\pi_{0}}=\int_{\varphi_{1}}\omega_{\xi},
\]
we remark that $s\mapsto q$ gives the mirror map for the family $W$
of elliptic curves. Similarly, if we set 
\[
\ell(Q):=\mathcal{T}:=\frac{1}{(2\pi\sqrt{-1})^{3}}\int_{\mathcal{M}(3\varphi_{0})}\eta,
\]
then $s\mapsto Q$ is the local mirror map for $Y$. The initial VMHS
$\mathcal{V}$ is that on $H^{3}(Y)$, with integral basis%
\footnote{We will ignore for now the fact that $\gamma_{1}$ is really $\frac{1}{3}$
of an integral class; it is a more convenient choice for our purposes
than $\mathcal{M}(\varphi_{1})^{\vee}$.%
} $\gamma=\{\gamma_{3},\gamma_{2},\gamma_{1}\}$ where 
\[
\gamma_{3}:=\mathbb{T}^{\vee},\;\gamma_{2}=\mathcal{M}(3\varphi_{0})^{\vee},\;\gamma_{1}=\mathcal{M}(3\varphi_{1})^{\vee}.
\]
From the exact sequence we can read off the weight filtration
\[
W_{6}=\mathcal{V}\supset W_{5}=W_{4}=W_{3}=\langle\gamma_{2},\gamma_{1}\rangle=\text{im}\{\mu\}\supset W_{2}=\{0\},
\]
and Hodge filtration (except for $\mathcal{F}^{3}=\langle\eta\rangle$).
The extension data are recorded by $\mathcal{T}=\langle AJ(\Xi),3\varphi_{0}\rangle$
and $\Phi:=\langle AJ(\Xi),3\varphi_{1}\rangle$.

The monodromy logarithm 
\[
[N]_{\gamma}=\left(\begin{array}{ccc}
0 & 0 & 0\\
-1 & 0 & 0\\
\frac{1}{2} & -1 & 0
\end{array}\right)
\]
leads to a relative weight filtration $M_{\bullet}$. The resulting
$\mathcal{V}_{rel}$ has Hodge-Tate basis \begin{flalign*}
\; & e_3 := \frac{\eta}{(2\pi\sqrt{-1})^3}= \gamma_3+\mathcal{T}\gamma_2 + \Phi \gamma_1 \in \mathcal{F}^3\cap M_6, &\\
\; & e_2 :=  \mu(\omega) = \gamma_2 +\tau \gamma_1 \in \mathcal{F}^2 \cap M_4, &\\
\; & e_1 = \gamma_1 \in \mathcal{F}^1\cap M_2. 
\end{flalign*}From transversality 
\[
\gamma_{2}+\frac{d\Phi}{d\mathcal{T}}\gamma_{1}=\nabla_{\partial_{\mathcal{T}}}e_{3}\in\mathcal{F}^{2}
\]
we deduce that $\frac{d\Phi}{d\mathcal{T}}=\tau$, which may also
be derived from the fact that logarithmic derivatives of the extension
classes give periods%
\footnote{That is, we have $\delta_{s}\mathcal{T}=\frac{1}{2\pi i}\pi_{0}$,
$\delta_{s}\Phi=\frac{1}{2\pi i}\pi_{1}$.%
} of $\tilde{\omega}_{\xi}$ {[}op. cit.{]}:
\[
\frac{d\Phi}{d\mathcal{T}}=\frac{s\cdot d\Phi/ds}{s\cdot d\mathcal{T}/ds}=\frac{\pi_{1}}{\pi_{0}}=\tau.
\]
This equality has the important consequence 
\[
\Phi'':=\frac{d^{2}\Phi}{d\mathcal{T}^{2}}=\frac{d\tau}{d\mathcal{T}}=\frac{\delta_{s}(\pi_{1}/\pi_{0})}{\delta_{s}\mathcal{T}}=\frac{2\pi\sqrt{-1}(\pi_{0}\delta_{s}\pi_{1}-\pi_{1}\delta_{s}\pi_{0})}{\pi_{0}^{3}}=\frac{\mathcal{Y}}{\pi_{0}^{3}},
\]
where $\mathcal{Y}$ is the (suitably normalized) Yukawa coupling
for the family $\{W_{\xi}\}$ of elliptic curves. Noting as well that
$\nabla_{\partial_{\mathcal{T}}}e_{2}=\frac{d\tau}{d\mathcal{T}}e_{1}$,
we conclude that
\[
[\nabla]_{e}=d+\left(\begin{array}{ccc}
0 & 0 & 0\\
1 & 0 & 0\\
0 & \Phi'' & 0
\end{array}\right)\otimes\frac{dQ}{(2\pi\sqrt{-1})Q}
\]
where $e=\{e_{3},e_{2},e_{1}\}$.

Turning to the A-model, we shall seek a quantum interpretation of
$\nabla$. Before doing so, we remark that by \cite{Ho} and \cite{DK},
under the local mirror map $\Phi$ may be identified as the local
Gromov-Witten prepotential\begin{equation}\label{eqn: local pre}\Phi\equiv\frac{1}{2}\mathcal{T}^{2}-\frac{1}{(2\pi \sqrt{-1})^2}\sum_{d}3dN_{d}Q^{d}.
\end{equation}modulo lower order terms in $\mathcal{T}$.%
\footnote{A different form of this result is already present in $\S$6.2 of
\cite{CKYZ}, about which we shall say more in the next section.%
} Differentiating \eqref{eqn: local pre} twice, we have
\[
1-\sum_{d}3d^{3}N_{d}Q^{d}=\frac{\mathcal{Y}}{\pi_{0}^{3}},
\]
in which the right-hand side has a pole where the family $W$ degenerates.
Directly computing $\langle AJ(\Xi),\varphi_{0}\rangle$ at this singular
elliptic curve gives $\Im(\mathcal{T}_{0})=\frac{27\sqrt{3}}{8\pi^{2}}L(\chi_{-3},2)$
\cite{DK}, and hence $Q_{0}=|e^{2\pi\sqrt{-1}\mathcal{T}_{0}}|=e^{-2\pi\Im(\mathcal{T}_{0})}$
for the radius of convergence. This ties the asymptotic growth rate
\[
\limsup_{d\to\infty}|N_{d}|^{\frac{1}{d}}=e^{2\pi\Im(\mathcal{T}_{0})}
\]
of the local Gromov-Witten numbers directly to the Beilinson regulator
of an algebraic cycle.

For the quantum interpretation, we consider the dual VMHS $\mathcal{V}^{\vee}$
on $H_{3}(Y)$ under the pairing $H^{3}(Y)\times H_{3}(Y)\to H_{0}(Y)=\mathbb{Z}.$
The dual integral (flat) basis is of course 
\[
\gamma_{1}^{\vee}=\mathbb{T}^{3},\;\;\gamma_{2}^{\vee}=\mathcal{M}(3\varphi_{0}),\;\;\gamma_{1}^{\vee}=\mathcal{M}(3\varphi_{1}),
\]
and in the dual Hodge basis $e^{\vee}=\{e_{3}^{\vee},e_{2}^{\vee},e_{1}^{\vee}\}$
we have\begin{equation}\label{eqn local nabla}[\nabla]_{e^{\vee}}=d-\left(\begin{array}{ccc} 0& 1&0\\ 0&0  & \Phi''\\ 0&0&0\\ \end{array}\right)\otimes\frac{dQ}{2\pi\sqrt{-1}Q}.
\end{equation}Now recalling that $Y^{\circ}=K_{\mathbb{P}^{2}}$, Hosono \cite{Ho}
proposed a homological mirror map
\[
\text{mir}:\, K_{0}^{c}(Y^{\circ})\to H_{3}(Y,\mathbb{Z})
\]
from coherent sheaves with compact support to homology classes of
Lagrangian 3-cycles, given explicitly by \begin{equation}\label{hosono maps}\mathcal{O}_{p}\mapsto\gamma_{3}^{\vee},\;\;\mathcal{O}_{\mathbb{P}^{1}}(-1)\mapsto\gamma_{2}^{\vee},\;\;\mathcal{O}_{\PP^{2}}(-2)\mapsto\gamma_{1}^{\vee}.
\end{equation}(The sheaves are all supported on the zero-section $\mathbb{P}^{2}\subset Y^{\circ}$.)
Making the identifications $e_{3}^{\vee}=[p]$, $e_{2}^{\vee}=[\mathbb{P}^{1}]$,
$e_{1}^{\vee}=[\mathbb{P}^{2}]$ under $\overline{\text{mir}}:\, H_{\text{even}}(Y^{\circ})\overset{\cong}{\to}H_{3}(Y)$,
we impose as before an integral structure on the A-model side by means
of the quantum deformed classes
\[
([p]=)\,[p]_{\mathcal{Q}}:=\gamma_{3}^{\vee},\;\;[\mathbb{P}^{1}]_{\mathcal{Q}}:=\gamma_{2}^{\vee},\;\;[\mathbb{P}^{2}]_{\mathcal{Q}}:=\gamma_{1}^{\vee}.
\]
Together with the filtrations $W_{-6}=W_{-5}=W_{-4}=\langle[p]\rangle\subset W_{-3}=H_{\text{even}}$,
and $\langle[p]\rangle=\mathcal{F}^{-3}\cap M_{-6}$, $\langle[\mathbb{P}^{1}]\rangle=\mathcal{F}^{-2}\cap M_{-4}$,
$\langle[\mathbb{P}^{2}]\rangle=\mathcal{F}^{-1}\cap M_{-2}$, this
determines the A-model (relative) variation matching that on the B-model.

Finally, consider the formal quantum product\begin{equation}\label{eqn: product}\begin{matrix}
e_1^{\vee}*e_3^{\vee}=0,\; e_1^{\vee}*e_2^{\vee}=-3e_3^{\vee},\; e_1^{\vee}*e_1^{\vee}=-3\Phi'' e_2^{\vee}, \\
e_2^{\vee}*e_3^{\vee}=0,\; e_3^{\vee}*e_3^{\vee}=0,\; e_2^{\vee}*e_2^{\vee}=0,
\end{matrix}
\end{equation}where we continue to identify classes under $\overline{\text{mir}}$.
This is compatible with the ordinary cup product in the sense that
\begin{flalign*} 
\; & e_1^{\vee}\cup e_3^{\vee}=[\mathbb{P}^2]\cup [p]=0, &\\
\; & e_1^{\vee}\cup e_2^{\vee} = [\mathbb{P}^2]\cup[\mathbb{P}^1]=(\mathbb{P}^2\cdot\mathbb{P}^1)_{Y^{\circ}}[p]=-3[p]=-3e_3^{\vee},\;\text{and} &\\
\; & e_1^{\vee}\cup e_1^{\vee} = [\mathbb{P}^2]\cup[\mathbb{P}^2]=-3[\mathbb{P}^1]=-3e_2^{\vee}, 
\end{flalign*}the last of which contains the leading term of $-3\Phi''=-3+\cdots$.
\begin{prop}
\label{prop: nabla}With the product \eqref{eqn: product}, \eqref{eqn local nabla}
may be rewritten
\[
\nabla=d+\left(\frac{1}{3}e_{1}^{\vee}*\right)\otimes\frac{dQ}{2\pi\sqrt{-1}Q}
\]
in terms of the quantum product with the zero-section $\mathbb{P}^{2}\subset K_{\mathbb{P}^{2}}$.
\end{prop}
This motivates the following
\begin{problem}
\label{LQC problem}Develop a general theory of quantum cohomology
for the local setting that produces $\nabla$ on $H_{\text{even}}(Y^{\circ})$
as Prop. \ref{prop: nabla}.
\end{problem}
We will obtain a solution for our running example in the next section.

The Abel-Jacobi maps from \cite{DK} touched on above may be viewed
as maps from $K_{2}^{\text{alg}}(W)=K_{2}(Coh(W))$ to ($\mathbb{C}/\mathbb{Z}(2)$-valued)
functionals on (classes of) Lagrangian 1-cycles on $W$. Noting that
$W^{\circ}$ is also an elliptic curve, we propose
\begin{problem}
\label{prob B}Derive (in general) a homological mirror to $AJ$.
This would produce a ``symplectic regulator'' map from $K_{2}(Fuk(W^{\circ}))$
to functionals on coherent sheaves on $W^{\circ}$. The functional
mirroring the $AJ$ class in our example would send $\mathcal{O}_{p}\mapsto\frac{(2\pi\sqrt{-1})^{2}}{3}\mathcal{T}$
and $\mathcal{O}_{W^{\circ}}\mapsto\frac{(2\pi\sqrt{-1})^{2}}{3}\Phi$.
\end{problem}
The motivation for such a quantum $AJ$ map is clear: it would bring
Beilinson's conjectures directly to bear upon the arithmetic of GW
invariants, in the context of the A-model VHS on quantum cohomology.
A first step might be to construct, in our example, a mirror in $K_{2}(Fuk(W^{\circ}))$
to the toric symbol $\{x,y\}\in K_{2}^{\text{alg}}(W)$ (i.e. the
higher cycle), by representing $K_{2}^{\text{alg}}(W)$ using the
Quillen category of $Coh(W)$ and applying homological mirror symmetry
for elliptic curves.

\section{Closed to Local}

We begin by summarizing a computation from \cite{CKYZ}. The setting
is a 2-parameter family $X_{\xi_{1},\xi_{2}}$ of $h^{2,1}=2$ CY
3-folds over a product of punctured disks, with $\hat{\eta}\in\Omega^{3}(X)$.
The mirror ($h^{1,1}=2$) CY has an elliptic fibration 
\[
X^{\circ}\overset{\bar{\rho}}{\to}\mathbb{P}^{2}
\]
with
\begin{itemize}
\item zero-section $D_{2}\cong\mathbb{P}^{2}$,
\item a line $C_{2}\cong\mathbb{P}^{1}\subset D_{2}$ with preimage $D_{1}=\bar{\rho}^{-1}(C_{2})$,
and
\item a fiber $C_{1}=\bar{\rho}^{-1}(p)$.
\end{itemize}
We will use the bases
\[
\left\{ \begin{array}{ccc}
J_{1}=[D_{2}]+3[D_{1}],\; J_{2}=[D_{1}] & \text{for} & H^{1,1}(X^{\circ})\\
C_{1},\; C_{2} & \text{for} & H^{2,2}(X^{\circ})
\end{array}\right.
\]
which are dual under cup product. The period vector for $\hat{\eta}$
takes the form
\[
\left(\Pi_{0},\tau_{1}\Pi_{0},\tau_{2}\Pi_{0},\partial_{\tau_{1}}\tilde{\Phi},\partial_{\tau_{2}}\tilde{\Phi},2\tilde{\Phi}-\delta_{\tau_{1}}\tilde{\Phi}-\delta_{\tau_{2}}\tilde{\Phi}\right)
\]
where $\Pi_{0}$ is the ``holomorphic period'' and\begin{equation*}\small
\xymatrix{\tilde{\Phi}:=\frac{3}{2}\tau_{1}^{3}+\frac{3}{2}\tau_{1}^{2}\tau_{2}+\frac{1}{2}\tau_{1}\tau_{2}^{2}+\left\{ \frac{17}{4}\tau_{1}+\frac{3}{2}\tau_{2}+C\right\} +\frac{1}{(2\pi\sqrt{-1})^{3}}\sum_{d_{1},d_{2}}\tilde{N}_{d_{1},d_{2}}q_{1}^{d_{1}}q_{2}^{d_{2}}
}
\normalsize
\end{equation*}is the prepotential.%
\footnote{This is the usual G-W prepotential \emph{plus} the bracketed lower-order
correction terms.%
} Here, $q_{j}=e^{2\pi\sqrt{-1}\tau_{j}}$ are the disk-coordinates
and $\tilde{N}_{d_{1},d_{2}}$ is the G-W invariant of the class $d_{1}[C_{1}]+d_{2}[C_{2}]$
on $X^{\circ}$; the K\"ahler class is simply $\tau_{1}J_{1}+\tau_{2}J_{2}$.

Now we take $\tau_{1}\to i\infty$ ($q_{1}\to0$) considered as the
``large volume limit'' for the fibers of $\bar{\rho}$. For the
purposes of G-W theory on the A-model, in this limit $X^{\circ}$
is equivalent to the total space of $\mathcal{N}_{D_{2}/X^{\circ}}\cong K_{\mathbb{P}^{2}}$,
i.e. $Y^{\circ}$ in the last section (with the map $\rho:Y^{\circ}\twoheadrightarrow\PP^{2}$).
On the B-model, which we shall henceforth ignore, the periods remaining
finite are $\Pi_{0}$, $\tau_{2}\Pi_{0}$, and\begin{equation}\label{eqn: finite}
(\partial_{\tau_1}-3\partial_{\tau_2})\tilde{\Phi}\; = \; \frac{1}{2}\tau_2^2-\frac{1}{4}+\frac{1}{(2\pi\sqrt{-1})^2}\sum_{d_1,d_2}\tilde{N}_{d_1,d_2}(d_1-3d_2)q_1^{d_1}q_2^{d_2}.
\end{equation}Indeed, actually taking the limit of \eqref{eqn: finite} (and writing
$\mathcal{T}:=\tau_{2}$, $Q:=e^{2\pi\sqrt{-1}\mathcal{T}}$, $N_{d}:=\tilde{N}_{0,d}$)
defines the local prepotential 
\[
\Phi_{\text{loc}}:=\frac{1}{2}\mathcal{T}^{2}-\frac{1}{4}-\frac{1}{(2\pi\sqrt{-1})^{2}}\sum_{d}3dN_{d}Q^{d}
\]
in agreement with \eqref{eqn: local pre}.%
\footnote{In fact, by a computation in \cite{Ho}, $\Phi=\Phi_{\text{loc}}-\frac{1}{2}\mathcal{T}+\frac{1}{2}$.%
}

The next step is to consider the limit of the quantum products of
classes in $H^{\text{even}}(X^{\circ})$ which come from $H_{c}^{\text{even}}(Y^{\circ})$($\cong H_{\text{even}}(Y^{\circ})$),
namely $[p]$, $[C_{2}]$, and 
\[
[D_{2}]=J_{1}-3J_{2}.
\]
In general, the only interesting products (not given by the cup product)
are 
\[
J_{j}*J_{k}=\sum_{\ell}\left(\partial_{\tau_{j}}\partial_{\tau_{k}}\partial_{\tau_{\ell}}\tilde{\Phi}\right)[C_{\ell}].
\]
So (using \eqref{eqn: finite}) we have 
\[
[D_{2}]*[D_{2}]=\left(\partial_{\tau_{1}}-3\partial_{\tau_{2}}\right)^{2}\left(\partial_{\tau_{1}}\tilde{\Phi}[C_{1}]+\partial_{\tau_{2}}\tilde{\Phi}[C_{2}]\right)
\]
\[
=-3[C_{2}]+\sum_{d_{1},d_{2}}\tilde{N}_{d_{1},d_{2}}(d_{1}-3d_{2})^{2}(d_{1}[C_{1}]+d_{2}[C_{2}])q_{1}^{d_{1}}q_{2}^{d_{2}},
\]
whereupon taking the limit $\lim_{q_{1}\to0}[D_{2}]*[D_{2}]=$
\[
\left\{ -3+\sum_{d}N_{d}(-3d)^{2}dQ^{d}\right\} [C_{2}]=
\]
\[
-3\left\{ 1-\sum_{d}3d^{3}N_{d}Q^{d}\right\} [C_{2}]
\]
gives $[\mathbb{P}^{2}]*[\mathbb{P}^{2}]=-3\Phi''[\mathbb{P}^{1}]$,
which is exactly what we wanted.

This makes a case for the general principle that the ``local restriction''
of the quantum product in a closed CY should remain finite under an
appropriate large volume limit. Beyond establishing this, a solution
to Problem \ref{LQC problem} would have to show the result is consistent
with a formula of the shape%
\footnote{For $Y^{\circ}\cong K_{\PP}\to\PP$ with $\PP$ a toic Fano surface,
negativity of $K_{\PP}$ allows us to express the local invariants
as closed invariants $\langle\alpha,\beta,\phi^{k}\rangle_{0,3,\iota_{*}(d)}$for
$\overline{Y^{\circ}}:=\PP(\mathcal{O}\oplus K_{\PP})\overset{\iota}{\supset}Y^{\circ}$,
cf. \cite[sec. 9]{CI}.%
}\begin{equation}\label{local * product}\alpha*_{\text{loc}}\beta:=\sum_{k}\sum_{d}\langle\alpha,\beta,\phi^{k}\rangle_{0,3,d}^{\text{loc}}\phi_{k}e^{\langle d,\mathcal{T}\rangle}
\end{equation}for $\alpha,\beta\in H_{even}\cong H_{c}^{even}$, $\mathcal{T}\in H^{2}$,
$d\in H_{2}$, and $\phi^{k}$ resp. $\phi_{k}$ dual bases of $H^{even}$
resp. $H_{even}.$

The resulting local quantum cohomology would then provide a direct
A-model approach to ``most'' of the variation of mixed Hodge structure
(the $\{I^{p,q}\}$ and $\nabla$-flat structure), leaving only the
\begin{problem}
\label{local quantum Z-str}Extend Iritani's construction of an integral
structure on $\nabla$-flat sections to the local CY setting.
\end{problem}
This is easily accomplished in our running example by ``taking LMHS
along $q_{1}=0$'' of the $\ZZ$-VHS over $(\Delta^{*})^{2}$ (common
to both the A- and B-models). More precisely, if $T_{1}$ denotes
the monodromy about $q_{1}=0$, with logarithm $N_{1}$, then the
limiting variation of MHS takes the form\[\includegraphics[scale=0.5]{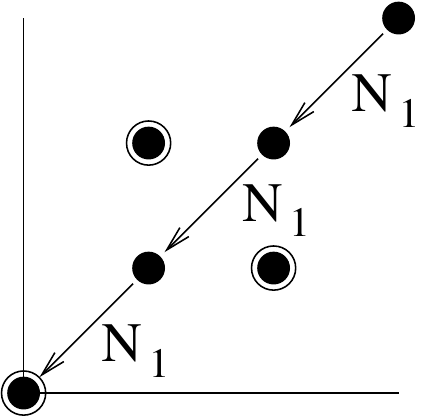}\]where
the circled bullets denote $\ker(N_{1})=\ker(T_{1}-\text{id})$. For
our purposes, then it will suffice to compute the limit of the $T_{1}$-invariant
``cycles'' in the $\hat{\Gamma}$-integral structure on the closed
A-model VHS $H^{even}(X^{\circ})$. 

Indeed, together with the Clemens-Schmid sequence, the assumption
that ``$Y^{\circ}$ is the A-model limit of $X^{\circ}$'' implies
that \begin{equation*}\tiny \xymatrix{
0 \to H_3(Y)(-2) \ar@{=} [d]  \ar [r] &  \underset{q_{1}\to0}{\lim}H^{3}(X)(1) \ar@{=} [d]\ar [r]^{N_1} &  \underset{q_{1}\to0}{\lim}H^{3}(X) \ar [r] \ar@{=} [d]& H^3(Y)  \ar@{=} [d] \to 0 \\
0 \to  H_{even}(Y^{\circ})(-2)   \ar [r] &  \underset{q_{1}\to0}{\lim}H^{even}(X^{\circ})(1) \ar [r]^{N_1} &  \underset{q_{1}\to0}{\lim}H^{even}(X^{\circ}) \ar [r] & H^{even}(Y^{\circ}) \to 0
}
\normalsize
\end{equation*}is an exact sequence of VMHS (in $q_{2}$). Iritani's procedure necessarily
gives integral $\nabla$-flat sections $\{\hat{\gamma_{i}}\}_{i=1}^{6}$
in $H^{even}(X^{\circ})$, with $\{\hat{\gamma}_{4},\hat{\gamma}_{5},\hat{\gamma}_{6}\}\subset\text{im}(N_{1})$
and $\{\hat{\gamma}_{1}^{\vee},\hat{\gamma}_{2}^{\vee},\hat{\gamma}_{3}^{\vee}\}\subset\ker(N_{1})$,
such that 
\[
\frac{\hat{\eta}}{\Pi_{0}}=\hat{\gamma}_{1}+\tau_{2}\hat{\gamma}_{2}+\left\{ (\partial_{\tau_{1}}-3\partial_{\tau_{2}})\tilde{\Phi}\right\} \hat{\gamma}_{3}+\sum_{j=4}^{6}\hat{\pi}_{j}(\underline{\tau})\hat{\gamma_{j}}.
\]
Taking the limit whilst killing $\text{im}(N_{1})$, then making the
change of basis $\{\hat{\gamma}_{1},\hat{\gamma}_{2},\hat{\gamma}_{3}\}=:\{\gamma_{1}+\frac{1}{2}\gamma_{3},\gamma_{2}-\frac{1}{2}\gamma_{3},\gamma_{3}\}$,
recovers
\[
e_{3}=\gamma_{1}+\mathcal{T}\gamma_{2}+\Phi\gamma_{3}
\]
in $H^{even}(Y^{\circ})$.

Of course, in analogy to \eqref{local * product}, it would be better
to solve Problem \ref{local quantum Z-str} in a manner intrinsic
to the local A-model. That is, there should be a direct construction
as in \eqref{gamma Z-str} assigning flat sections of $H_{even}(Y^{\circ})$
to classes in $K_{0}^{c}(Y^{\circ})$, and ``compatible with monodromy''.
In our example, \eqref{hosono maps} has this compatibility, since
$\otimes\mathcal{O}_{Y^{\circ}}(-J_{2})$ on the coherent sheaves
and monodromy about $q=0$ on the cycles have the same matrix
\[
\left(\begin{array}{ccc}
1 & 1 & 0\\
0 & 1 & 1\\
0 & 0 & 1
\end{array}\right).
\]

Apparently, either solution still leaves us a long way from the ``holy
grail'' of Problem \ref{prob B}.

\section{Open string}

Problem \ref{prob B} is probably intractable without major theoretical
developments. However, its rough analogue in the \emph{relative} situation
studied by Morrison and Walcher \cite{MW} appears to be more accessible.
In particular, there is nothing mysterious about the mirror of the
(usual, not higher) algebraic cycle -- it is just a Lagrangian.

The B-model in the example we consider (following {[}op. cit.{]})
comprises:
\begin{itemize}
\item $X=$ a double-cover of the mirror quintic family, with holomorphic
form $\omega\in\Omega^{3}(X)$;
\item $Z\in CH^{2}(X)_{\text{hom}}$ a family of algebraic 1-cycles (for
analogy to $\S2$, think ``$K_{0}(Coh(X))$''); and
\item $\langle AJ_{X}^{2}(Z),\omega\rangle=$ the resulting ``truncated
normal function'', solving
\item the inhomogeneous Picard-Fuchs equation $D_{\text{PF}}^{\omega}\langle AJ_{X}^{2}(Z),\omega\rangle=:g$.
\end{itemize}
On the A-model side these data mirror to:
\begin{itemize}
\item $X^{\circ}$= the Fermat quintic;
\item $Z^{\circ}\cong\mathbb{RP}^{3}$ the real quintic, viewed as a Lagrangian
3-cycle (think ``$K_{0}(Fuk(X))$''); and
\item the Gromov-Witten generating function whose coefficients count holomorphic
disks bounding on $Z^{\circ}$,
\end{itemize}
which (under the mirror map) solves the same PF equation. 

As in the closed and local stories, GW numbers are therefore obtained
as power-series coefficients of a Hodge-theoretic function, with (in
this latter role) the Yukawa coupling replaced by the truncated normal
function.
\begin{problem}
Work out (in analogy with $\S\S$1-2) the $[\nabla]_{e}$ story. This
will require the \emph{full} normal function (not considered in {[}op.
cit.{]}), which means computing also $\langle AJ_{X}^{2}(Z),\nabla_{\partial_{\tau}}\omega^{3,0}\rangle$. 
\end{problem}
Since the B-model VMHS is an extension of the constant variation $\mathbb{Z}(-2)$
by the pure VHS $H^{3}(X)$, the extension class is defined over $\mathbb{R}$
hence given completely by $\langle AJ_{X}^{2}(Z),\omega^{3,0}\rangle$
and $\langle AJ_{X}^{2}(Z),\nabla_{\partial_{\tau}}\omega^{3,0}\rangle$.
The extension arises geometrically from the residue exact sequence\begin{equation*}
\small
\xymatrix{0 \ar [r] & H^3(X) \ar [r] \ar @{=} [d] & H^3(X\setminus |Z|) \ar [r] & \ker \left( H^4_{|Z|}(X) \to H^4(X)\right) \ar [r] & 0 \\
0 \ar [r] & H^3(X) \ar [r] & \mathbb{E} \ar @{^(->} [u] \ar [r] & \mathbb{Q}(-2) \ar @{^(->} [u] \ar [r] & 0}
\normalsize
\end{equation*}Completely missing, however, is an approach to the following.
\begin{problem}
Can one produce the extension class from the pair $X^{\circ}$,$Z^{\circ}$
from the standpoint of quantum cohomology and the A-model VMHS?
\end{problem}
To illustrate its difficulty, a naive attempt to mirror the exact
sequence approach, viz.
\[
0\to\frac{H^{\text{even}}(X^{\circ})}{H_{Z^{\circ}}^{6}(X^{\circ})}\to H^{\text{even}}(X^{\circ}\backslash Z^{\circ})\to\ker\left(H_{Z^{\circ}}^{3}(X^{\circ})\to H^{3}(X^{\circ})\right)\to0,
\]
fails due to the vanishing of the third term. The result of {[}op.
cit.{]}, however, that the ``truncated'' extension class is given
by the open GW generating function, gives one reason to believe the
problem has interesting content.
\begin{rem}
We briefly note another interesting phenomenon that arises in the
open setting, related to Remark \ref{rem arithmetic}. Even on a family
of CY 3-folds defined oveer $\QQ$, algebraic cycles often force an
algebraic extension $L/\QQ$ upon us, as in the case of the van Geemen
lines on the mirror quintic family studied by Laporte and Walcher
\cite{LW}. The resulting limits of truncated normal functions can
then often be expressed in terms of the Borel regulator on $K_{3}^{ind}(L)$
(see \cite{GGK2} for the theoretical reason). This makes the open
setting ideal terrain for exploring generalizations of the A-model
$\hat{\Gamma}$-construction where the B-model LMHS does not correspond
to a $\QQ$-rational limiting motive.
\end{rem}

\subsection{Local to open}

Recent work of Chan, Lau, Leung, Tseng and Wu \cite{CLL,CLT,LLW}
has brought to light an interesting relation between the (local) mirror
map and certain open Gromov-Witten invariants for a toric Calabi-Yau
manifold $Y^{\circ}$. The first three authors conjecture in \cite{CLL}
that the SYZ mirror construction (applied to $Y^{\circ}$) inverts
the mirror map given by a normalized integral basis of single-log-divergent
periods of the Hori-Vafa mirror $Y$. With the integrality hypothesis
dropped, the conjecture is established in \cite{CLT} for $Y^{\circ}=K_{Z}$
with $Z$ a compact toric Fano variety; it is known integrally for
toric surfaces \cite{LLW} and a handful of other examples \cite{CLL},
including $K_{\PP^{2}}$.

We briefly describe the case $Y^{\circ}=K_{\PP^{2}}$ in the notation
of $\S2$. Take $\beta_{0}$ to denote the class of a holomorphic
disk bounding on the zero section $D\,(\cong\PP^{2})\subset Y^{\circ}$,
$\ell$ the class of a line $L\,(\cong\PP^{1})\subset D$; and let
$\mathcal{T}[\rho^{-1}(L)]\in H^{2}(Y,\CC)$ be th K\"ahler class
with corresponding K\"ahler parameter $Q=e^{2\pi i\mathcal{T}}$.
Then the SYZ construction in {[}op. cit{]} produces produces the noncompact
Calabi-Yau in $(\CC^{*})^{2}\times\CC^{2}$ given by \begin{equation}\label{eqn 5.1}UV=c(Q)+X+Y+\frac{Q}{XY},
\end{equation}where $c(Q)=1+\sum_{k\geq1}n_{\beta_{0}+k\ell}Q^{k}$ is a local Gromov-Witten
generating series. An easy change of coordinates exhibits \eqref{eqn 5.1}
as the Hori-Vafa manifold $Y_{\xi}$ of \eqref{eqn HV mir}, with
$\xi=-\frac{Q}{c(Q)^{3}}$; taking the cube gives \begin{equation}\label{eqn 5.2}s(Q)=-\frac{Q}{c(Q)^3}.
\end{equation}The observation of {[}op. cit.{]} is that \eqref{eqn 5.2} inverts
the local mirror map
\[
Q(s)=e^{2\pi\sqrt{-1}\mathcal{T}(s)}=\exp\left(\frac{1}{(2\pi\sqrt{-1})^{2}}\int_{\mathcal{M}(3\varphi_{0})}\eta\right)
\]
in $\S2$. So just as for $\Phi$, we have an enumerative interpretation
for $\mathcal{T}$, and one can use the computation%
\footnote{up to the sign and term $\frac{1}{2}$ which are required for consistency
with $\S2$ and \cite{Ho}%
}
\[
\mathcal{T}(s)=\ell(s)+\frac{1}{2}+\frac{1}{2\pi\sqrt{-1}}\sum_{k\geq1}\frac{\binom{3k}{k,k,k}}{k}s^{k}
\]
in \cite{CLL} or \cite{DK} to compute $c(Q)=1-2Q+5Q^{2}-32Q^{3}+\cdots.$

We conclude with one final
\begin{problem}
Can one use the formulae in $\S5$ of \cite{DK} for the integral
periods of Hori-Vafa mirrors, to establish integrality in \cite{CLT}? \end{problem}

\curraddr{${}$\\
\noun{Department of Mathematical and Statistical Sciences}\\
\noun{University of Alberta, Canada}}

\email{${}$\\
\emph{e-mail}: doran@math.ualberta.ca}

\curraddr{\noun{${}$}\\
\noun{Department of Mathematics, Campus Box 1146}\\
\noun{Washington University in St. Louis}\\
\noun{St. Louis, MO} \noun{63130, USA}}

\email{\emph{${}$}\\
\emph{e-mail}: matkerr@math.wustl.edu}
\end{document}